\documentclass[11pt]{article}

\usepackage{geometry}   
\usepackage{graphicx}
\usepackage{amsmath}
\usepackage{amssymb}
\usepackage{authblk}
\usepackage{amsthm}
\usepackage{epstopdf}
\usepackage{enumitem}
\usepackage{lscape}
\usepackage{array}
\usepackage{adjustbox}
\usepackage{multirow}
\usepackage{bm}
\usepackage{tikz}
\usetikzlibrary{matrix, arrows, decorations.pathreplacing, shapes}
\setlist[description]{leftmargin=\parindent,labelindent=\parindent}

\numberwithin{equation}{section}

\newtheorem{thm}{Theorem}[section]

\theoremstyle{remark}

\theoremstyle{definition}
\newtheorem{defn}[thm]{Definition}

\newcommand{\Pf}{\operatorname{Pf}}
\newcommand{\Spec}{\operatorname{Spec}}
\newcommand{\Proj}{\operatorname{Proj}}

\newcommand{\exc}{\operatorname{exc}}

\newcommand{\CC}{\mathbb{C}}

\newcommand{\PP}{\mathbb{P}}
\newcommand{\QQ}{\mathbb{Q}}

\newcommand{\frakm}{\mathfrak{m}}

\title{Divisorial extractions from singular curves in a \\ smooth 3-fold, II: low codimension.}
\author{Tom Ducat}
\affil{RIMS, Kyoto University, Kyoto 606-8502, Japan. \\ \tt{taducat@kurims.kyoto-u.ac.jp}}
\date{\today}     
                       
\begin{document}

\maketitle

\begin{abstract}
Following the first paper, we continue to study Mori extractions from singular curves centred in a smooth 3-fold. We treat the case where the divisorial extraction exists in relative codimension at most 3.
\end{abstract}

\section*{Introduction}
\addcontentsline{toc}{section}{Introduction}

The focus of this paper is on the classification of Mori extractions---that is 3-dimensional divisorial extractions over $\CC$ with at worst terminal singularities. The first paper \cite{duc1} introduced a method to construct the coordinate ring of a divisorial extraction from $C\subset X$, a singular curve $C$ contained in a smooth 3-fold $X$, using type I Gorenstein unprojection. This divisorial extraction is the unique Mori extraction from $C$, if such an extraction exists, and is isomorphic to the blowup of the \emph{symbolic power algebra} 
\[ \sigma\colon Y \cong \Proj_X \bigoplus I_{C/X}^{[n]} \to X\]
where $I_{C/X}$ is the ideal defining $C\subset X$ (see \cite{duc1} Proposition 1.4). Thus, if $\bigoplus I_{C/X}^{[n]}$ is generated by $r+1$ generators, we get a presentation of $Y$ as a codimension $r$ subvariety $Y\subset X\times w\PP^r$ for some weighted projective space $w\PP^r$ of dimension $r$. 

This paper studies the cases for which a Mori extraction exists in codimension $\leq3$. It may be unreasonable to expect an explicit classification in the completely general case. In fact we know that the divisorial extractions can have arbitrarily high codimension -- see for example \cite{thesis} \S6.3 where one such family of examples is constructed by serial unprojection. Moreover, in general it seems hard to find an exact condition on $C\subset X$ for which the divisorial extraction has at worst terminal singularities.

Nevertheless, the hope is that this list will find useful applications in studying the explicit birational geometry of 3-folds, e.g.\ in the Sarkisov program. For instance, Prokhorov \& Reid \cite{prokreid} used a Sarkisov link beginning with the simplest Mori extraction in Table \ref{non_lci_table} (type $\bm{A}^1_1$) to construct a new $\QQ$-Fano 3-fold of index 2, and more complicated divisorial extractions can be used to construct other interesting examples \cite{duc3}.

\subsection{Main result}
The aim of the paper is to prove the following result:
\begin{thm}
Suppose $P\in C\subset S_X\subset X$ is an inclusion of varieties, as in \S1, and suppose $\sigma\colon Y\subset X\times w\PP^r\to X$ is a codimension $r$ Mori extraction from $C\subset X$ constructed by unprojection. Then: 
\begin{enumerate}
\item If $r\leq 2$ then $C\subset X$ is a l.c.i.\ and $Y$ is the ordinary blowup of $C$. The possible numerical types for $C\subset S_X$ are given in Table \ref{lci_table}, along with the curves extracted by $\sigma|_{S_Y}$. Moreover, every case in the table has a Mori extraction.
\item If $r=3$ then $w\PP=\PP(1,1,1,2)$ and $Y$ has a singularity of index 2. The possible numerical types for $C\subset S_X$ are given in Table \ref{non_lci_table}, along with the curves extracted by $\sigma|_{S_Y}$. 
\end{enumerate}
\end{thm}

(We do not claim that every case appearing in Table \ref{non_lci_table} has a Mori extraction.)

\subsection{Notation and conventions}
\label{notation}

\paragraph{Resolutions of Du Val singularities.} We fix the following numbering of the $ADE$ Dynkin diagrams:
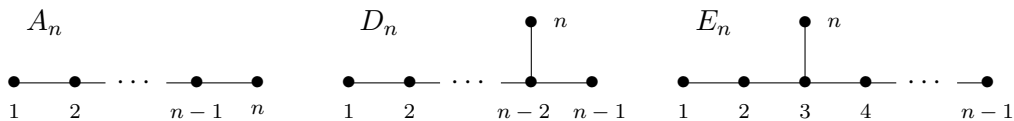
\begin{center} \begin{tikzpicture}[scale=0.8]
	\node at (-1.5,1) {$A_n$};
	\draw (-2,0) -- (-0.5,0)
	 	 (0.5,0) -- (2,0);
	\node[label={[label distance=-2]below:\scriptsize 1}] at (-2,0) {$\bullet$};
	\node[label={[label distance=-2]below:\scriptsize 2}] at (-1,0) {$\bullet$};
	\node at (0,0) {$\cdots$};
	\node[label={[label distance=-2]below:\scriptsize $n-1$}] at (1,0) {$\bullet$};
	\node[label={[label distance=-2]below:\scriptsize $n$}] at (2,0) {$\bullet$};
	
	\node at (4,1) {$D_n$};
	\draw 
		(3.5,0) -- (5,0)
		(6,0) -- (7.5,0)
		(6.5,0) -- (6.5,1);
	\node[label={[label distance=-2]below:\scriptsize 1}] at (3.5,0) {$\bullet$};
	\node[label={[label distance=-2]below:\scriptsize 2}] at (4.5,0) {$\bullet$};
	\node at (5.5,0) {$\cdots$};
	\node[label={[xshift=-0.1cm,label distance=-2]below:\scriptsize $n-2$}] at (6.5,0) {$\bullet$};
	\node[label={[xshift=0.1cm,label distance=-2]below:\scriptsize $n-1$}] at (7.5,0) {$\bullet$};
	\node[label={[label distance=-2]right:\scriptsize $n$}] at (6.5,1) {$\bullet$};
	
	\node at (9.5,1) {$E_n$};
	\draw 
		(9,0) -- (12.5,0)
		(13.5,0) -- (14,0)
		(11,0) -- (11,1);
	\node[label={[label distance=-2]right:\scriptsize $n$}] at (11,1) {$\bullet$};
	\node[label={[label distance=-2]below:\scriptsize 1}] at (9,0) {$\bullet$};
	\node[label={[label distance=-2]below:\scriptsize 2}] at (10,0) {$\bullet$};
	\node[label={[label distance=-2]below:\scriptsize 3}] at (11,0) {$\bullet$};
	\node[label={[label distance=-2]below:\scriptsize 4}] at (12,0) {$\bullet$};
	\node at (13,0) {$\cdots$};
	\node[label={[label distance=-2]below:\scriptsize $n-1$}] at (14,0) {$\bullet$};
	
\end{tikzpicture} \end{center}
Given the minimal resolution of a Du Val singularity $\mu\colon (E\subset \widetilde{S}) \to (P\in S)$ we write $E_i$ for the exceptional divisor corresponding to the $i$th vertex in the Dynkin diagram, $\widetilde{C}_i\subset \widetilde{S}$ for a smooth curve transverse to $E_i$ and $C_i=\mu(\widetilde{C}_i)\subset S$.

\paragraph{Numerical types.} Given $C\subset S$, the \emph{numerical type of $C$} is $C\equiv_{\text{num}} \sum_{i=1}^n \alpha_iC_i$ where $\widetilde{C}\equiv_{\text{num}} \sum_{i=1}^n \alpha_i\widetilde{C}_i$ on $\widetilde{S}$. (Note that $\widetilde{C}$ may intersect $E$ nontransversely.)


\paragraph{The cycle $\Delta_i$.}
When $P\in S$ is of type $D_n$ we define the numerical type $\Delta_i\subset S$, for $i=1,\ldots, n$, by the following formula:
\[ \Delta_i = \begin{cases} C_i & 1\leq i \leq n-2 \\ C_{n-1} + C_{n} & i=n-1 \\ 2C_{n-1} \text{ or } 2C_{n} & i=n \end{cases} \]
At one point (see Table \ref{non_lci_table}) we need to distinguish between the two different cases for $\Delta_n$, but otherwise we treat both cases simultaneously.

\subsection{Acknowledgements.}

The author is a International Research Fellow of the Japanese Society for the Promotion of Science and this work was supported by Grant-in-Aid for JSPS Fellows, No.\ 15F15771. 

He would also like to thank Miles Reid for reading a preliminary draft of this paper and for his help in improving it.

\section{The divisorial extraction}

As in \cite{duc1}, we assume the following inclusion of algebraic varieties over $\CC$:
\[ P\in C \subset S_X \subset X \]
where $P\in X$ is a smooth point in a 3-fold, $P\in S_X$ is a Du Val singularity and $C$ is a curve with a singularity at $P$. We are interested in the existence of a Mori extraction $\sigma \colon (F \subset Y)\to (C\subset X)$, where $F$ is the exceptional divisor extracted from $C$. By \cite{duc1} Proposition 1.4, $Y$ is isomorphic to the blowup of the \emph{symbolic power algebra} of the ideal $I_{C/X}$ defining $C\subset X$.

\paragraph{The general elephant.} The existence of the Du Val hypersurface $S_X$ is a consequence of the \emph{general elephant conjecture} which states that, for a Mori extraction (or a Mori flipping contraction) $\sigma\colon Y\to X$, a general member $S_Y\in|{-K_Y}|$ and $S_X:=\sigma(S_Y)\in|{-K_X}|$ should both have at worst Du Val singularities. Then, by adjunction, the restriction $\sigma|_{S_Y}\colon S_Y\to S_X$ is a partial crepant resolution.

Unlike \cite{duc1}, we \emph{do not} assume that $S_X$ is a general hypersurface section containing $C$. Part of the information given in the description of a Mori extraction is which curves (if any) are extracted from $P\in S_X$ by $\sigma|_{S_Y}$.

\subsection{Constructing the divisorial extraction}

We briefly recall the method explained in \cite{duc1} \S2.4 to construct $Y$.

\paragraph{A normal form for $C$.}
If $C\subset S_X$ is not a local complete intersection (l.c.i.) then by \cite{duc1} Lemma 2.1 we can write the equations of $C$ as the minors of a $2\times 3$ matrix $M$:
\begin{equation}
\bigwedge^2 \begin{pmatrix}  \:\:\phi\: & \begin{matrix} -b \\ a \end{matrix}  \end{pmatrix} = 0 
\label{C_eqns}
\end{equation}
where $\phi$ is a $2\times2$ matrix such that $\det \phi$ is the equation of the Du Val singularity $P\in S_X$ and $a,b
\in\CC[x,y,z]$ are some functions on $X$. After a change of variables we can take $\phi$ to be one of the following matrices:
\[ \begin{array}{ccc}
\bm{A}_{n-1}^{i} \colon \begin{pmatrix} x & y^i \\ y^{n-i} & z \end{pmatrix} \quad & 
\bm{D}^l_{n} \colon \begin{pmatrix} x & y^2 + z^{n-2} \\ z & x \end{pmatrix} \quad & 
\bm{D}^r_{2k} \colon \begin{pmatrix} x & yz + z^k \\ y & x \end{pmatrix} \\ \\
\bm{D}^r_{2k+1} \colon\begin{pmatrix} x & yz \\ y & x + z^k \end{pmatrix} \quad & 
\bm{E}_6 \colon\begin{pmatrix} x & y^2 \\ y & x + z^2 \end{pmatrix} \quad & 
\bm{E}_7\colon \begin{pmatrix} x & y^2 + z^3 \\ y & x \end{pmatrix} 
\end{array} \]
($E_8$ does not appear in this list as it is factorial and has no such nontrival matrix factorisation.) We can use row and column operations to ensure that $b=b(y,z)$ and $a=a(y,z)$ (or $a=a(x,y)$ in case $\bm{A}_{n-1}^i$).

We write $I := I(M)$ for the ideal generated by the entries of $M$ and $I(\phi)$ for the ideal generated by those of $\phi$. Clearly $I(\phi)\subseteq I$ and it is also clear that $C\subset X$ is a l.c.i. if (and only if) $I\not\subseteq \frakm_P$ (i.e.\ at least one of $a,b\notin\frakm_P$).

\paragraph{The blowup of $C$.}

We can write down the (ordinary) blowup of $C\subset X$ as a complete intersection of codimension 2, with equations given by the two syzygies coming from Cramer's rule:
\begin{equation} 
\sigma_0\colon Y_0\subset X\times \PP^2_{(\nu:\xi:\eta)}\to X \quad\quad \begin{pmatrix}  \:\:\phi\: & \begin{matrix} -b \\ a \end{matrix}  \end{pmatrix} \begin{pmatrix}  \nu \\ -\xi \\ \eta  \end{pmatrix} = 0 
\label{Y_eqns}
\end{equation}
In particular $\eta$ corresponds to the equation $\det \phi$ defining $S_X$. If $I\subseteq \frakm_P$ the exceptional divisor of this blowup has two components $F_0\cup D_0$, where $F_0$ is a reduced divisor dominating $C$ and $D_0=V(I)$ is a (not necessarily reduced) codimension 1 subscheme dominating $P\in X$. 

\begin{defn}
The \emph{width of $D_0$} is defined to be $w(D_0) := \dim_\CC \CC[x,y,z] / I$. 
\end{defn}

Since $I(\phi)\subseteq I$, the width of $D_0$ is bounded by $w_\phi := \dim_\CC \CC[x,y,z] / I(\phi)$. In the different cases this bound is given by:
\[ \renewcommand{\arraystretch}{1.3} 
\begin{array}{c|cccccc}
\text{Type} & \bm{A}_{n-1}^{i} & \bm{D}^l_{n} & \bm{D}^r_{2k} & \bm{D}^r_{2k+1} & \bm{E}_6 & \bm{E}_7 \\ \hline
I(\phi) & (x,y^i,z) & (x,y^2,z) & (x,y,z^k) & (x,y,z^k) & (x,y,z^2) & (x,y,z^3) \\
w_\phi & i & 2 & k & k & 2 & 3
\end{array} \]

In particular $D_0\cong \PP^2\times \Spec \CC[\varepsilon]/(\varepsilon^{w(D_0)})$ is the largest possible subscheme contained in $Y_0$ supported on $P\times \PP^2$ (c.f.\ \cite{pagoda} Definition 5.3).

\paragraph{Unprojecting $D_0$.} Since $Y_0$ and $D_0$ are both complete intersections in $X\times \PP^2$ (and hence Gorenstein) we can use the unprojection theorem (\cite{par} Theorem 1.5) to unproject $D_0\subset Y_0$. In this simple case we can obtain the unprojection variable $\zeta$ by rewriting \eqref{Y_eqns} as a $2\times 3$ matrix annihilating $I$ and applying \cite{kino} Example 4.16. By construction $\zeta$ satisfies $\zeta I \subseteq I_{C/X}^2$ so $\zeta$ is naturally an element of the second symbolic power $\zeta\in I_{C/X}^{[2]}$. Hence unprojection gives us a Gorenstein variety $Y_1\subset X\times \PP(1,1,1,2)$ in codimension 3 with a birational map
\[ u \colon (D_0\subset Y_0) \dashrightarrow (Q\in Y_1) \]
$Y_1$ is defined by five equations and so, by a theorem of Buchsbaum and Eisenbud, we can (and will) write the equations as the maximal Pfaffians of a $5\times 5$ skew matrix. This is also the way the equations are presented in \cite{kino} Example 4.1. 

Unprojection is a type of birational surgery that blows up $D_0$ to a Cartier divisor $\widetilde{D_0}$ and contracts the birational transform. In the nicest case $Y_0$ has only isolated ordinary nodal singularities at the points where $D_0$ fails to be Cartier and blowing up $D_0$ makes a small resolution of these nodes. However unprojection can be more complicated in practice. For example, $Y_0$ may contain a line of ordinary nodes (or a line of transverse $A_n$ singularities) along $D_0$, and then unprojection introduces a new divisor $D_1\subset Y_1$ above $D_0$. This can lead to a chain of \emph{serial unprojections}. One family of divisorial extractions given an arbitrarily long sequence of serial unprojections is constructed in \cite{thesis} \S6.3.

\section{Codimension $\leq2$ cases}

We first treat the cases in which a Mori extraction $\sigma\colon Y\to X$ exists in codim $\leq 2$. In all of these cases $C\subset X$ is a l.c.i.\ and then, by a theorem of Cutkosky \cite{cut} (see also \cite{duc1}, Lemma 1.6), it is known that a necessary and sufficient condition for $C$ to have a Mori extraction is that $C$ is contained in a smooth hypersurface. Using this condition we give a classification of the possible numerical types for $C$ and we also classify which curve (if any) is extracted from $S_X$.

\subsection{Summary}

We summarise all the cases in Table \ref{lci_table}. The table lists the $ADE$ type of $P\in S_X$, the format for $C\subset S_X$, which exceptional curves (if any) are extracted by $\sigma|_{S_Y}\colon S_Y\to S_X$ and the possible numerical types of $C$. 

Each case appears in the table up to a symmetry of the corresponding Dynkin diagram, e.g.\ for the $\bm{E}_6$ case the divisorial extraction from $C_5$ (with $\exc(\sigma|_{S_Y})=E_1$) also includes the divisorial extraction from $C_1$ (with $\exc(\sigma|_{S_Y})=E_5$).

Since they are contained in a smooth hypersurface section, all the cases in the table have a Mori extraction, even the degenerate cases in which $\widetilde{C}\cap E$ is nontransverse. 

\begin{table}[htp]
\renewcommand{\arraystretch}{1.3}
\caption{Curves with codimension $\leq2$ Mori extraction.}
\begin{center}
\begin{tabular}{|c|c|c|m{7.8cm}|} \hline
Type & Format & $\exc(\sigma|_{S_Y})$ & Numerical types \\ \hline
\multirow{2}{*}{$A_{n-1}$} & l.c.i. & $\emptyset$ & $\sum_{j=1}^{n-1} \alpha_jC_j$ : $\sum_j j\alpha_j = n$ or  $\sum_j (n-j)\alpha_j = n$ \\
& $\bm{A}_{n-1}^i$ & $E_{n-i}$ & $\sum_{j=1}^{n-1} \alpha_jC_j$ : $\sum_j j\alpha_j = i$ or  $\sum_j (n-j)\alpha_j = n-i$ \\ \hline
\multirow{3}{*}{$D_n$} & l.c.i. & $\emptyset$ & $\Delta_{2i}$ $(\forall i\leq\lfloor\tfrac{n}{2}\rfloor)$, $C_1+\Delta_{2i-1}$ $(\forall i\leq\lfloor\tfrac{n+1}{2}\rfloor)$ \\
& $\bm{D}_n^l$ & $E_1$ & $\Delta_{2i-1}$ $(\forall i\leq\lfloor\tfrac{n+1}{2}\rfloor)$ \\ 
& $\bm{D}_n^r$ & $E_n$ & $C_n, \: C_1+C_{n-1}$ \\ \hline
\multirow{2}{*}{$E_6$} & l.c.i. & $\emptyset$ & $C_6, \: C_1+C_5, \: C_3, \: C_1+C_2, \: 3C_1, \: C_4+C_5, \: 3C_5$ \\ 
& $\bm{E}_6$ & $E_1$ & $C_5, \: C_2, \: 2C_1$ \\ \hline
\multirow{2}{*}{$E_7$} & l.c.i. & $\emptyset$ & $C_1, \: C_5, \: 2C_6, \: C_2, \: C_6+C_7$ \\ 
& $\bm{E}_7$ & $E_6$ & $C_6, \: C_7$ \\ \hline
$E_8$ & l.c.i. & $\emptyset$ & $C_7, \: C_1, \: C_6, \: C_8$ \\ \hline
\end{tabular} 
\smallskip

Note: The cycles $\Delta_j$ appearing in the $D_n$ cases are defined in \S0.2.
\end{center}
\label{lci_table}
\end{table}%

\subsection{Proof of the classification}

The proof follows from explicit calculations. First note that one of the following occurs:
\begin{itemize}
\item[(i)] $C\subset S_X$ is a l.c.i., in which case $S_Y\cong S_X$.
\end{itemize}
We write $C=V(f,g)\subset X$ where $S_X=V(f)$ and $g$ is the equation of a smooth hypersurface. The Mori extraction from $C$ is given by the codimension 1 model: 
\[ Y=V(f\xi - g\eta)\subset X\times \PP^1_{(\xi:\eta)}. \]
$S_Y=V(\eta)\subset Y$ meets the central fibre $Z = P\times \PP^1$ at the point $P_\xi$, where all variables apart from $\xi$ vanish, and $P_\xi \in S_Y$ is a singularity with equation $f=0$. So $\sigma|_{S_Y}$ is an isomorphism $S_Y \cong S_X$. 

We can now classify the possible numerical types for $C$ by writing down a minimal resolution $\mu\colon (E\subset \widetilde{S}_X)\to (P\in S_X)$ and calculating $\widetilde{C}\cap E$, for all $C$ subject to the condition that $g\notin \frakm_P^2$. (See \S\ref{exampleCalc} for an example of the kind of necessary calculation.)
\begin{itemize}
\item[(ii)] $C\subset S_X$ is not a l.c.i., in which case $S_Y\not\cong S_X$.
\end{itemize}
We use \eqref{C_eqns} to write the equations of $C\subset X$ and \eqref{Y_eqns} to write the equations of $Y=Y_0$. Of the two terms $a,b$ appearing in the format \eqref{Y_eqns} we must have at least one of $a,b\notin\frakm_P$, else $Y$ contains an unprojection divisor. The central fibre of $\sigma$ is $Z=P\times \PP^1_{(\xi:\nu)}$ and $Z\subset S_Y=V(\eta)$. Hence $\exc(\sigma|_{S_Y})\neq \emptyset$ and $S_Y\not\cong S_X$.

We can now check in each of the cases which curve is extracted from $P\in S_X$ and which numerical types are possible for $C$ by explicitly calculating $\widetilde{C}\subset\widetilde{S}_X$, subject to condition that at least one of $a,b\notin \frakm_P$.

\subsubsection{$A_{n-1}$ cases}
\label{type_A}

In the $A_{n-1}$ case we can consider $S_X$ as the $\tfrac1n(1,n-1)$ cyclic quotient singularity
\[ \pi \colon \CC^2_{u,v} \to \CC^2 / \bm{\mu}_n =: S_X \]
where $\bm{\mu}_n = \langle \varepsilon : \varepsilon^n = 1 \rangle$, the cyclic group of the $n$th roots of unity, acts on $\CC^2$ by $(u,v) \mapsto (\varepsilon u, \varepsilon^{n-1} v)$. We write $x,y,z=u^n,uv,v^n$ for the invariants of this action which satisfy the relation $xz = y^r$. We can pull back any curve $C\subset S_X$ to an invariant curve $\Gamma := \pi^{-1}(C)\subset \CC^2_{u,v}$ given by a semi-invariant \emph{orbifold equation} $\Gamma=V\big(\gamma(u,v)\big)$. If $C\equiv\sum_{j=1}^{n-1}\alpha_jC_j$ then $\gamma$ factors (analytically) as a product
\[ \gamma(u,v) = \prod_{j=1}^{n-1}\gamma_j(u^j,v^{n-j}) \]
where $\gamma_j(U,V)\in\CC[[x,y,z]][U,V]$ is a homogeneous polynomial of degree $\alpha_j$ whose roots correspond to the intersection points of $\widetilde C\cap E_j$ counted with multiplicity.
 
By the normal form \eqref{C_eqns}, the equations of $C\subset X$ are given by the following format for some $i$:
\[ \bigwedge^2\begin{pmatrix}
x & y^i & -b \\
y^{n-i} & z & a
\end{pmatrix} =0 \]
In this case the orbifold equation is given by $\gamma(u,v) = au^i + bv^{n-i}$ and the equations $ax+by^{n-i}$, $ay^i+bz$ are given by rendering the invariants $u^{n-i}\gamma$, $v^i\gamma$ in terms of $x,y,z$. 

If $a_0,b_0$ are the constant terms of $a,b$, then the condition that $C$ is contained in a smooth hypersurface is equivalent to asking that at least one of $a_0,b_0\neq0$. But $a_0$ is the coefficient of $u^i$ in $\gamma$ so, by the product expression above, $a_0\neq0$ if and only if $i = \sum_{j=1}^{n-1} j\deg{\gamma_j} = \sum_{j=1}^{n-1} j\alpha_j$. Similarly $b_0\neq0$ if and only if $n-i = \sum_{j=1}^{n-1} (n-j)\alpha_j$.

\subsubsection{Example calculation}
\label{exampleCalc}

We explain the $\bm{E}_6$ case as an example of the calculations that can be used to verify the other type $D$ and $E$ cases appearing Table \ref{lci_table}.

\paragraph{The equations of $C$.}
Using the normal form \eqref{C_eqns} in the $\bm{E}_6$ case we write the equations of $C\subset X$ as the $2\times 3$ minors of the matrix:
\[ \bigwedge^2\begin{pmatrix}
x & y^2 & -b \\
y & x + z^2 & a
\end{pmatrix} \]
for some choice of $a,b\in\CC[y,z]$. As before we write $a_0,b_0$ for the constant terms of $a,b$. Note that $C$ is contained in a smooth hypersurface if and only if at least one of $a_0,b_0\neq0$. 

\paragraph{Explicit resolution of $P\in S_X$}
We resolve $S_X=V\big(x(x+z^2) - y^3\big)$ as follows. Let $b_y$, $b_z$ be the two following coordinate changes:
\[ b_y \colon (x,y,z) \mapsto (xy,y,yz), \quad b_z \colon (x,y,z) \mapsto (xz,yz,z). \]
Then the minimal resolution $\mu\colon (E\subset \widetilde{S}_X) \to (P\in S_X)$ can be covered by the five following affine charts:
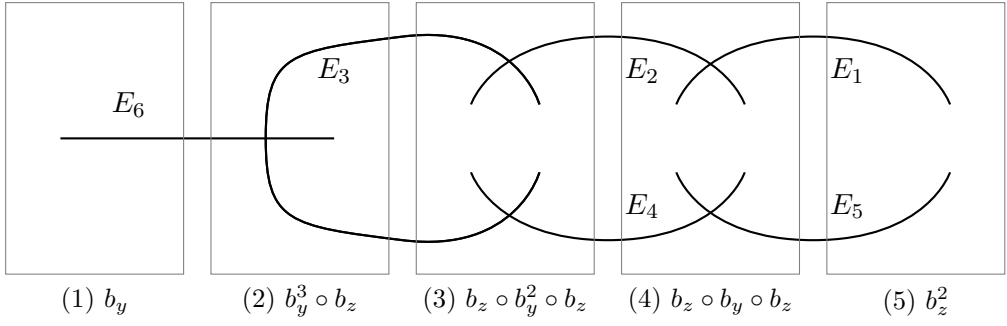
\begin{center} \begin{tikzpicture}[scale=1.8]
  	\draw[thick] (0,0) -- (2,0);
	\draw[thick] plot [smooth, tension=1.5] coordinates {(3.5,-0.25) (2.5,-0.75) (1.5,0) (2.5,0.75) (3.5,0.25)};
	\draw[thick] plot [smooth, tension=1.5] coordinates {(3.5,-0.25) (2.5,-0.75) (1.5,0) (2.5,0.75) (3.5,0.25)};
	\draw[thick] plot [smooth, tension=1.5] coordinates {(3,0.25) (4,0.75) (5,0.25)};
	\draw[thick] plot [smooth, tension=1.5] coordinates {(4.5,0.25) (5.5,0.75) (6.5,0.25)};
	\draw[thick] plot [smooth, tension=1.5] coordinates {(3,-0.25) (4,-0.75) (5,-0.25)};
	\draw[thick] plot [smooth, tension=1.5] coordinates {(4.5,-0.25) (5.5,-0.75) (6.5,-0.25)};
	
  	\draw[gray] (-0.4,1) rectangle (0.9,-1);
  	\draw[gray] (1.1,1) rectangle (2.4,-1);
  	\draw[gray] (2.6,1) rectangle (3.9,-1);
  	\draw[gray] (4.1,1) rectangle (5.4,-1);
  	\draw[gray] (5.6,1) rectangle (6.9,-1);
	
	\node at (5.75,0.5) {$E_1$};
	\node at (4.25,0.5) {$E_2$};
	\node at (2,0.5) {$E_3$};
	\node at (4.25,-0.5) {$E_4$};
	\node at (5.75,-0.5) {$E_5$};
	\node at (0.5,0.25) {$E_6$};
	
	\node at (0.25,-1.2) {\small (1) $b_y$};
	\node at (1.75,-1.2) {\small (2) $b_y^3\circ b_z$};
	\node at (3.25,-1.2) {\small (3) $b_z\circ b_y^2\circ b_z$};
	\node at (4.75,-1.2) {\small (4) $b_z\circ b_y\circ b_z$};
	\node at (6.25,-1.2) {\small (5) $b_z^2$};
\end{tikzpicture} \end{center}
For example, we reach the chart (4) by the change of coordinates 
\[ \mu_4 \colon (x,y,z) \mapsto (xy^2z^4,y^2z^3,yz^2). \]
In this chart $\widetilde{S}_X=V\big (x(x+1) - y^2z \big)$ is a smooth surface and (up to a choice of relabelling) the exceptional curves are given by
\[ E_1 = V(x,y), \quad E_2 = V(x,z), \quad E_4 = V(x+1,z), \quad E_5 = V(x+1,y). \]

\paragraph{Description of $C\subset S_X$.}
By calculating $\widetilde{C}\subset \widetilde{S}_X$ we see that there are exactly three numerical types for curves satisfying the condition that one of $a_0,b_0\neq0$. We do the computation in chart (4) of the resolution given above and leave the rest of the calculation to the reader. In fact, in each of these three cases $\widetilde{C}\cap E=\emptyset$ outside of chart (4) (apart from possibly the point at $\infty$ on either $E_1$ or $E_5$ contained in chart (5)).

\begin{description}
\item[(i)] First assume $b_0\neq 0$. Then in chart (4) $\widetilde{C}$ is given by
\[ \widetilde{C}=V\big(a'xz + b', \: x(x+1) - y^2z) \] 
(where $a' = \mu_4(a)$ and $b'=\mu_4(b)$). Since $b_0\neq 0$, $\widetilde{C}\cap E_i=\emptyset$ for $i=1,2,4$ and $\widetilde{C}$ intersects $E_5$ at one point transversely at the point $x=y=a_0z - b_0=0$ (if $a_0=0$ then this intersection point is the point at $\infty$ in chart (5)). Hence $C\equiv C_5$.

\item[(ii)] If $b_0=0$ then we replace $b$ by $by+cz$ and now we assume $a_0,c_0\neq 0$. Then in chart (4) $\widetilde{C}$ is given by
\[ \widetilde{C}=V \big( a'y+b'(x+1)yz + c'(x+1), \: a'x+b'y^2z^2+c'yz \big) \] 
$\widetilde{C}$ intersects $E_2$ transversely at the point $x=z=a_0y+c_0=0$. Moreover,  if $c_0\neq0$ then $\widetilde{C}$ does not intersect $E_1$ and if $a_0\neq0$ then $\widetilde{C}$ does not intersect $E_i$ for $i=3,4,5,6$. Hence $C\equiv C_2$. 

\item[(iii)] Lastly we assume that $c_0=0$ and $a_0\neq 0$ in case (ii), so we replace $by+cz$ by $by+cz^2$. Then in chart (4) $\widetilde{C}$ is given by
\[ \widetilde{C}=V \big( a'+b'(x+1)z + c'(x+1)z^2, \: x(x+1) - y^2z \big) \] 
$\widetilde{C}$ intersects $E_1$ twice according to the roots of $c_0z^2+b_0z+a_0=0$ (again, if $c_0,b_0=0$ then these points may be the point at $\infty $ in chart (5)). Moreover, if $a_0\neq0$ then $\widetilde{C}$ does not intersect any other component of $E$, so $C\equiv 2C_1$. 
\end{description}

In all other cases $C$ is not contained in a smooth hypersurface and $Y_0$ contains an unprojection plane.

\paragraph{The general elephant $S_Y$.} 
The last thing to check is that $\exc(\sigma|_{S_Y})=E_1$ as claimed. By \eqref{Y_eqns} we write $Y=Y_0\subset X\times \PP^2_{(\xi:\nu:\eta)}$ as a complete intersection: 
\[ Y = V\big(y^2\xi - x\nu + b\eta, \: y\nu - (x+z^2)\xi + a\eta\big) \]
with central fibre $Z=\PP^1_{\xi:\nu}\subset S_Y$. Therefore exactly one curve is extracted from $P\in S_X$.

It follows directly from these equations that $S_Y$ is smooth apart from at the point $P_\xi$, where all variables except $\xi$ vanish. At $P_\xi$ we can eliminate $x$ by the equation $x=y\nu - z^2$ to be left with the $D_5$ singularity 
\[ V(y^2 - y\nu^2 + \nu z^2)\subset \CC^3_{y,z,\nu}. \]
Hence $\sigma|_{S_Y}$ extracts either $E_5$ or $E_1$ from $P\in S_X$ and the extracted curve is independent of the choice of $a,b$. To see which one it is we can consider $C' = (\sigma|_{S_Y})^{-1}_*C$, the birational transform of $C$ under $\sigma|_{S_Y}$. For example take case (iii) above, we see that $C'$ is the complete intersection of $S_Y$ and $y=a+b\nu+c\nu^2$, so that $C'$ intersects $Z$ twice according to the two roots of $a_0+b_0\nu+c_0\nu^2=0$. Therefore it is $E_1$ that is extracted.

%

\section{Codimension 3 cases}

Now we treat the cases in which a Mori extraction $\sigma\colon Y \to X$ exists in codimension 3. As described in \S1.1, such an extraction has a model 
\[ \sigma \colon Y := Y_1 \subset X\times \PP(1,1,1,2) \to X \]
where $Y_1$ is the unprojection of $D_0\subset Y_0$. We consider each format $\bm{A}_{n-1}^{i},\ldots,\bm{E}_7$ in turn and split into subcases depending on $w(D_0)$, the width of the unprojection divisor $D_0\subset Y_0$ appearing in the construction of $Y$. 

Using the equations of $Y$ we give conditions for the reduced central fibre $Z= \sigma^{-1}(P)_{\text{red}}$ to be small (i.e.\ purely 1-dimensional) and conditions for $Y$ to have isolated singularities. Then we use these conditions to give a classification of the possible numerical types for $C\subset S_X$. We also check directly which curves, if any, are extracted by $\sigma|_{S_Y}$ from $P\in S_X$.

We are not claiming at present that every case in the list gives rise to a Mori extraction. A full analysis of the singularities appearing in each case would be very involved, and we would have to consider a lot of exceptional cases in which $\widetilde{C}$ degenerates to a curve intersecting $E$ nontransversely. However, once we establish that $\sigma|_{S_Y}\colon S_Y \to S_X$ is a partial crepant resolution of a Du Val singularity then the only non-terminal singularities of $Y$ must have centre in $Y\setminus S_Y$, so it would be enough just to check the singularities in this open set.

\subsection{Summary}

We summarise all the cases in Table \ref{non_lci_table}. In this table we list the format of $C\subset S_X$, the width $w=w(D_0)$ of the unprojection divisor $D_0\subset Y_0$, which exceptional curves (if any) are extracted by $\sigma|_{S_Y}\colon S_Y\to S_X$ and the possible numerical types of $C$. 

As before, each case appears in the table only up to a symmetry of the corresponding Dynkin diagram.

In the $\bm{A}_{n-1}^{i}$ case the table only contains the generic case for the numerical type of $C$. For a more detailed description of the possible numerical types see \S\ref{A_cases}.

\begin{table}[h]
\caption{Curves with codimension 3 Mori extraction.}
\begin{center}
\renewcommand{\arraystretch}{1.3}
\begin{tabular}{|c|c|c|m{5cm}|} \hline
Format & $w(D_0)$ & $\exc(\sigma|_{S_Y})$ & Numerical types \\ \hline
\multirow{3}{*}{$\bm{A}_{n-1}^i$}
& $w<i$ & $E_{n-i-w}+E_{n-i+w}$ & generically $C_w + C_i + C_{n-w}$ \\
& $i$ $(n>2i)$ & $E_{n-2i}$ & generically $2C_i + C_{n-i}$ \\ 
& $i$ $(n=2i)$ & $\emptyset$ & generically $3C_i$ \\ \hline
\multirow{2}{*}{$\bm{D}_n^l$}
& 1 & $E_2$ & $C_1 + \Delta_{2i}$ $(\forall i\leq \lfloor\tfrac{n}{2}\rfloor)$ \\
& 2 & $\emptyset$ & $2C_1 + \Delta_{2i-1}$ $(\forall i\leq\lfloor\tfrac{n+1}{2}\rfloor)$ \\ \hline
$\bm{D}_{n}^r$ & $w\leq \lfloor\tfrac{n}{2}\rfloor$ & $E_{n-2w}$ ($\emptyset$ if $n=2w$) & $C_n+\Delta_{2w}$ $(*)$, $C_1 + C_n + \Delta_{2w-1}$, $C_{n-1} + \Delta_{2w+1}$ $(**)$ \\ \hline
\multirow{2}{*}{$\bm{E}_6$}
& 1 & $E_2$ & $C_5+C_6$, $C_1+C_4$ \\
& 2 & $E_5$ & $C_1+2C_5$, $C_3+C_5$, $2C_4$ \\ \hline
\multirow{3}{*}{$\bm{E}_7$}
 & 1 & $E_5$ & $C_1+C_6$, $C_4$ \\
 & 2 & $E_1$ & $C_5+C_6$, $C_2+C_6$ \\
 & 3 & $\emptyset$ & $3C_6$, $2C_6+C_7$ \\ \hline
\end{tabular}
\begin{flushleft}
Note: There are two exceptions in the $\bm{D}_n^r$ case:
\begin{description}
\item[$(*)$] $C_n+\Delta_{2w}$ has $w(D_0)=w$ \emph{unless} $2w = n$, in which case $3C_n$ has $w(D_0)=w$ and $2C_{n-1}+C_n$ has $w(D_0)=w-1$. 
\item[$(**)$] $C_{n-1}+\Delta_{2w+1}$ has width $w(D_0)=w$ \emph{unless} $2w + 1 = n$, in which case $3C_{n-1}$ has $w(D_0)=w-1$ and $C_{n-1}+2C_n$ has $w(D_0)=w$. (If $w=2k$ then $2w+1>n$ and we ignore this case.)
\end{description}
\end{flushleft}
\end{center}
\label{non_lci_table}
\end{table}

\subsection{Proof of the classification}

We now prove the classification by dividing into cases according to the different formats $\bm{A}_{n-1}^{i},\ldots,\bm{E}_7$.

\subsubsection{Type $\bm{A}_{n-1}^i$}
\label{A_cases}

In this section we assume that our curve $C\subset S_X\subset X$ is of type $\bm{A}_{n-1}^i$. In this case $w(D_0)\leq i$ and we split into the following three subcases:
\begin{description}
\item[(i)] $w(D_0)<i$,
\item[(ii)] $w(D_0)=i$ and $2i<n$,
\item[(iii)] $w(D_0)=i$ and $2i=n$.
\end{description}

\paragraph{The equations of $C$ and $Y$.}
Our curve $C\subset X$ and the variety $Y\subset X\times \PP(1,1,1,2)$, given by the unprojection of $D_0\subset Y_0$, are defined by the following equations:
\[ \bigwedge^2\begin{pmatrix}
x & y^i & -(cy^w + dz) \\
y^{n-i} & z & ax + by^w 
\end{pmatrix} \quad \quad \Pf \begin{pmatrix}
\zeta & \nu & y^{i-w}\xi + c\eta & -d\eta \\
 & -a\eta & y^{n-i-w}\nu + b\eta & \xi \\
 & & z & y^w \\
 & & & x
\end{pmatrix} \]
where the five equations defining $Y$ are written as the Pfaffians of a skew matrix, as in \cite{kino} Example 4.1. (Only the strict upper diagonal part of the matrix is written). Without loss of generality, we collect terms together so that $a = a(x,y)$, $b = b(y)$, $c = c(y)$ and $d=d(y,z)$. Moreover we can assume that $n\geq 2i$ and that $w=w(D_0)$, so that if $w<i$ then at least one of $b_0,c_0\neq0$. Note that $i-w>0$ in case (i) and $i-w=0$ in cases (ii), (iii). Similarly, $n-i-w>0$ in cases (i), (ii) and $n-i-w=0$ in case (iii).

\paragraph{The general elephant $S_Y$.} 
\begin{description}
\item[(i)] $\exc(\sigma |_{S_Y}) = \PP^1_{(\xi:\zeta)} \cup \PP^1_{(\nu:\zeta)}$ corresponding to $E_{n-i+w}$ and $E_{n-i-w}$ respectively.
\end{description}
Setting $x,y,z,\eta=0$ in the equations defining $Y$ we see that $\exc(\sigma |_{S_Y})$ consists of the two irreducible components $\PP^1_{(\xi:\zeta)} \cup \PP^1_{(\nu:\zeta)}$. The restriction of $\xi,\nu,\zeta$ to $S_X$ can be written in terms of $u,v$ and $\gamma$ (the orbinates on $S_X$ introduced in \S\ref{type_A}) as follows: 
\[ \xi = [u^{n-i}\gamma], \quad \nu = [v^i\gamma], \quad \zeta = \frac{\xi\nu}{y^w} = [u^{n-i-w}v^{i-w}\gamma^2] \]
Then coordinates along $\PP^1(1,2)_{(\xi:\zeta)}$ are given by the ratio 
\[ (\xi^2:\zeta) = (u^{2(n-i)}\gamma^2:u^{n-i-w}v^{i-w}\gamma^2) = (u^{n-i+w}:v^{i-w}) \]
so this component corresponds to the exceptional divisor $E_{n-i+w}$ above $P\in S_X$. Similarly $\PP^1(1,2)_{(\nu:\zeta)}$ corresponds to $E_{n-i-w}$. Indeed, we can also use the equations to see that $S_Y$ is smooth apart from an $A_{i-w-1}$ singularity at $P_\xi$, an $A_{2w-1}$ singularity at $P_\zeta$ and an $A_{n-i-w-1}$ singularity at $P_\nu$, just as expected.

By a similar calculation in the other two cases:
\begin{description}
\item[(ii)] $\exc(\sigma |_{S_Y}) = \PP^1_{(\nu:\zeta)}$ corresponding to $E_{n-2i}$.
\item[(iii)] $\sigma |_{S_Y}$ is an isomorphism.
\end{description}

\paragraph{Small central fibre.} 
In all cases $Z$ is small and has $\leq3$ irreducible components, unless one of the following three conditions hold:
\begin{itemize}
\item $b_0=c_0=0$ in case (i): but this contradicts $w=w(D_0)$.
\item $a_0=b_0=0$ in either case (i) or (ii): then $Z$ contains a new unprojection divisor $D_1=V(x,y,z,\xi)$. 
\item $c_0=d_0=0$ in case (i) or $b_0=c_0=d_0=0$ in case (ii): then $Z$ contains a new unprojection divisor $D_1=V(x,y,z,\nu)$. 
\end{itemize}
In particular, in case (iii) $Z$ is always small.

\paragraph{Isolated singularities.}
In all of these cases the index 2 point $Q\in Y$ is a (hyper)quotient singularity:
\[ \big( \xi\nu = y^w + ad\eta^2 \big) \subset \CC^4_{\xi,\nu,\eta,y} \: / \: \tfrac12(1,1,1,0) \]

If $w=1$ then this is a $\tfrac12(1,1,1)$ singularity and $Y$ must have isolated singularities. A necessary condition for $Y$ to be terminal in this case is that not all of $a_0,b_0,c_0,d_0=0$, else $Z=\PP^1_{(\eta:\zeta)}$ and $Y$ has an index 1 singularity at $P_\eta$ with embedding dimension $\geq 5$, as in \cite{duc1} Lemma 2.3.

If $w>1$ then this is a $cA/2$ singularity, although possibly not isolated. Indeed if both $a_0,d_0=0$ then $Z$ contains $\PP^1_{(\eta:\zeta)}$ as a component and $Y$ is singular along this line.  

\paragraph{Description of $C\subset S_X$.} We summarise the necessary conditions for $Y$ to be the Mori extraction from $C$ in each of the three cases (i)--(iii).
\begin{description}
\item[(i)] At least one of $a_0,b_0\neq0$, at least one of $b_0,c_0\neq 0$, at least one of $c_0,d_0\neq 0$ and, if $w>1$, at least one of $a_0,d_0\neq 0$.
\item[(ii)] At least one of $a_0,b_0,c_0\neq 0$, at least one of $c_0,d_0\neq 0$ and, if $w>1$, at least one of $a_0,d_0\neq 0$.
\item[(iii)] At least one of $a_0,b_0,c_0,d_0\neq 0$ and, if $w>1$, at least one of $a_0,d_0\neq 0$.
\end{description}
In the case where $a,b,c,d$ are chosen generically then $C\subset S_X$ is the curve:
\[ \text{(i)} \:\: C_w + C_i + C_{n-w} \quad\quad \text{(ii)} \:\: 2C_i + C_{n-i} \quad\quad \text{(iii)} \:\: 3C_i \]
For example, in case (i) $C$ is given by the orbifold equation
\[ \gamma(u,v) = au^{n+i} + bu^{i+w}v^{w} + cu^{w}v^{n-i+w} + dv^{2n-i} \]
which, if all $a_0,b_0,c_0,d_0\neq0$, has initial term
\[ \gamma(u,v) = a_0(u^{n-w} + \tfrac{b_0}{a_0} v^{w})(u^i + \tfrac{c_0}{b_0} v^{n-i} )( u^{w} + \tfrac{d_0}{c_0} v^{n-w} ) + \cdots \]

\subsubsection{Type $\bm{D}_n^l$}

Now we assume that our curve $C\subset S_X\subset X$ is of type $\bm{D}_{n}^l$. In this case $w\leq 2$ so we split into the two subcases: (i) $w=1$, (ii) $w=2$.

\paragraph{The equations of $C$ and $Y$.}
Our curve $C\subset X$ and the variety $Y\subset X\times \PP(1,1,1,2)$ given by the unprojection of $D\subset Y_0$ are given by the following equations:
\[ \bigwedge^2 \begin{pmatrix}
x & y^2 + z^{n-2} & -(cy^w+dz) \\
z & x & ay^w+bz
\end{pmatrix} \quad\quad \Pf \begin{pmatrix}
\zeta & \nu & y^{2-w}\xi + c\eta & z^{n-3}\xi + d\eta \\ 
 & \xi & a\eta & \nu + b\eta \\
 & & -z & y^w \\
 & & & x
\end{pmatrix} \]
where we choose $w=w(D_0)$.

\paragraph{The general elephant $S_Y$.} 
In case (i) $\exc(\sigma |_{S_Y})=\PP^1_{(\xi:\zeta)}$ corresponding to $E_{2}$ and in case (ii) $\sigma |_{S_Y}$ is an isomorphism.

\paragraph{Small central fibre.}
In case (i) $Z$ is small unless either $a_0=c_0=0$ (which contradicts out choice of $w=w(D_0)$) or $c_0=d_0=0$, in which case $Y_1$ contains a new unprojection divisor $D_1=V(x,y,z,\nu)$. In case (ii) $Z$ is always small.

\paragraph{Isolated singularities.}
In case (i) $Q\in Y$ is a $\tfrac12(1,1,1)$ singularity and $Y$ must have isolated singularities. In case (ii) $Q\in Y$ is the hyperquotient singularity:
\[ \big( y^2 + \xi(z^{n-3}\xi + d\eta) = \nu(\nu + b\eta) \big) \subset \CC^4_{\xi,\nu,\eta,y} \: / \: \tfrac12(1,1,1,0) \]
If both $b_0,d_0=0$ then $\PP^1_{(\eta:\zeta)}\subseteq Z$ and $Y$ becomes singular along this line.  

\paragraph{Description of $C\subset S_X$.} Using the restrictions obtained above on the curves for which $Y$ is a Mori extraction, we can now give an explicit description of $C\subset S_X$.
\begin{description}
\item[(i)] At least one of $a_0,c_0\neq 0$ and at least one of $b_0,d_0\neq0$. Assume for the moment that $d_0\neq 0$. Resolving the $D_n$ singularity we see that $\widetilde{C}$ intersects $E_1\cong\PP^1_{(u_1:v_1)}$ according to $a_0u_1+c_0v_1=0$ and $E_2\cong\PP^1_{(u_2:v_2)}$ according to $c_0u_2+d_0v_2=0$. Moreover $\widetilde{C}$ does not meet any other exceptional curve. (Note that by our assumptions both these linear polynomials are nonzero and $\widetilde{C}$ meets the intersection point $E_1\cap E_2$ when $c_0=0$). Hence in this generic case $C\equiv C_1+C_2$, as represented in the following diagram (where the white node is the curve extracted from $S_X$): 
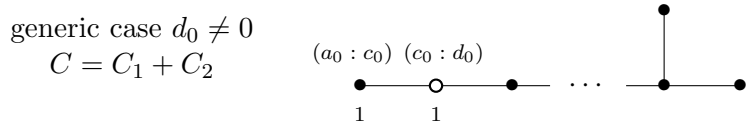
\begin{center} \begin{tikzpicture}[scale=1]
	\draw 
		(2.5,0) -- (5,0)
		(6,0) -- (7.5,0)
		(6.5,0) -- (6.5,1);
	\node[label={[label distance=-2]below:\scriptsize 1},label={[xshift=-0.1cm, label distance=-2]above:\scriptsize $(a_0:c_0)$}] at (2.5,0) {$\bullet$};
	\node[label={[label distance=-2]below:\scriptsize 1},label={[xshift=0.1cm, label distance=-2]above:\scriptsize $(c_0:d_0)$}] at (3.5,0) {$\bullet$};
	\node at (4.5,0) {$\bullet$};
	\node at (5.5,0) {$\cdots$};
	\node at (6.5,0) {$\bullet$};
	\node at (7.5,0) {$\bullet$};
	\node at (6.5,1) {$\bullet$};
	\draw[thick, fill=white] (3.5,0) circle (2.3pt) ; 
	\node at (-0.5,0.5) {$\begin{matrix} \text{generic case $d_0\neq 0$} \\ C=C_1+C_2 \end{matrix}$};
\end{tikzpicture} \end{center}
If $d_0=0$ then necessarily $c_0\neq0$, else $Z$ has a 2-dimensional component. We can check that making the replacement
\[ dz \mapsto \begin{cases}  
dz^i & i < \tfrac{n-1}{2} \\
dz^i & i = \tfrac{n-1}{2} \\
(dz \pm \sqrt{-1}c)z^{i-1} & i = \tfrac{n}{2}
\end{cases} \quad \text{gives} \quad C = \begin{cases} C_1 + C_{2i} \\ C_1 + C_{n-1} + C_n \\ C_1 + 2C_{n-1} \text{ or } C_1 + 2C_{n} \end{cases} \]
Moreover these are the only numerical types satisfying $c_0\neq0$. Using the cycle $\Delta_i$ defined in \S\ref{notation} we can write this family of curves as $C_1+\Delta_{2i}$, for $i\leq \tfrac{n}{2}$.
\item[(ii)] By a similar calculation we see that $C\equiv 3C_1$ if $d_0\neq 0$, intersecting $E_1\cong\PP^1_{(u_1:v_1)}$ according to the roots of the cubic $a_0u_1^3 + c_0u_1^2v_1 + b_0u_1v_1^2 + d_0v_1^3=0$. If $d_0=0$ then $b_0\neq0$ else $Y$ has non-isolated singularities. Then $C$ can degenerate to any of the curves $2C_1 + \Delta_{2i-1}$, for $i\leq \tfrac{n+1}{2}$.
\end{description}

\subsubsection{Type $\bm{D}_{n}^r$}

The cases $\bm{D}_{2k}^r$ and $\bm{D}_{2k+1}^r$ turn out to be very similar, even though the formats for $\phi$ in \eqref{C_eqns} initially look quite different. Therefore we only consider the even case $\bm{D}_{2k}^r$. Both cases are summarised in Table \ref{non_lci_table}.

We assume that our curve $C\subset S_X\subset X$ is of type $\bm{D}_{2k}^r$. In this case $w\leq k$ and we split into the two subcases: (i) $w<k$, (ii) $w=k$.

\paragraph{The equations of $C$ and $Y$.}
Our curve $C\subset X$ and the variety $Y\subset X\times \PP(1,1,1,2)$ given by the unprojection of $D\subset Y_0$ are given by the following equations:
\[ \bigwedge^2 \begin{pmatrix}
x & yz + z^k & -(cy + dz^w) \\
y & x & ay + bz^w
\end{pmatrix} \quad \quad \Pf \begin{pmatrix}
\zeta & \nu & z\xi + c\eta & z^{k-w}\xi + d\eta \\
 & \xi & \nu + a\eta & b\eta \\
 & & -z^w & y \\
 & & & x
\end{pmatrix} \]
where we choose $w=w(D_0)$.

\paragraph{The general elephant $S_Y$.} 
In case (i) $\exc(\sigma |_{S_Y})=\PP^1_{(\xi:\zeta)}$ corresponding to $E_{2(k-w)}$ and in case (ii) $\sigma |_{S_Y}$ is an isomorphism.

\paragraph{Small central fibre.}
In case (i) $Z$ is small unless either $b_0=d_0=0$ (which contradicts out choice of $w=w(D_0)$) or $c_0=d_0=0$, in which case $Y_1$ contains a new unprojection divisor $D_1=V(x,y,z,\nu)$. In case (ii) $Z$ is always small.

\paragraph{Isolated singularities.}
$Q\in Y$ is the hyperquotient singularity:
\[ \big( z^w + \nu(\nu + a\eta) =  \xi(z\xi + c\eta) \big) \subset \CC^4_{\xi,\nu,\eta,z} \: / \: \tfrac12(1,1,1,0) \]
If $w>1$ and both $a_0,c_0=0$ then $\PP^1_{(\eta:\zeta)}\subseteq Z$ and $Y$ becomes singular along this line.

\paragraph{Description of $C\subset S_X$.}
By explicitly resolving $S_X$ and computing $\widetilde{C}\subset\widetilde{S}_X$ we see that, for $w<k$, the three possible numerical equivalence classes for $C$ listed in Table \ref{non_lci_table} correspond to:
\begin{itemize}
\item $C_n + \Delta_{2w}$ : $c_0,d_0\neq0$ (the generic case) 
\item $C_1+C_n + \Delta_{2w-1}$ : $c_0=0$ and $a_0,d_0\neq0$
\item $C_{n-1} + \Delta_{2w+1}$ : $d_0=0$ and $b_0,c_0\neq0$
\end{itemize}
In all other cases either $c_0=d_0=0$ (so that $Z$ is not small), $b_0=d_0=0$ (so that $w(D_0)>w$) or $a_0=c_0=0$ (so that $Y$ has non-isolated singularities). Similarly for $w=k$ we have two cases to consider: the generic case $3C_n$, given by $c_0\neq0$, and a special case $C_1+C_{n-1}+2C_n$ given by $c_0=0$ and $a_0\neq0$.

\subsubsection{Type $\bm{E}_6$}

Now we assume that our curve $C\subset S_X\subset X$ is of type $\bm{E}_{6}$. In this case $w \leq 2$ and we split into the two subcases: (i) $w=1$, (ii) $w=2$.

\paragraph{The equations of $C$ and $Y$.}
Our curve $C\subset X$ and $Y\subset X\times \PP(1,1,1,2)$ are given by the following equations:
\[ \bigwedge^2 \begin{pmatrix}
x & y^2 & -(cy + dz^w) \\
y & x + z^2 & ay + bz^w
\end{pmatrix} \quad\quad  \Pf \begin{pmatrix}
\zeta & \nu & y\xi + c\eta & -d\eta \\
 & \xi & \nu + a\eta & z^{2-w}\xi - b\eta \\
 & & z^w & y \\ 
 & & & x
\end{pmatrix} \]
where we choose $w=w(D_0)$.

\paragraph{The general elephant $S_Y$.} In both cases $\exc(\sigma |_{S_Y})=\PP^1_{(\xi:\zeta)}$, corresponding to $E_2$ in case (i) and to $E_5$ in case (ii).

\paragraph{Small central fibre.}
In case (i) $Z$ is small unless either $b_0=d_0=0$ (contradicting $w=w(D_0)$) or $c_0=d_0=0$, in which case $Z$ contains a new unprojection divisor $D_1=V(x,y,z,\nu)$. In case (ii) $Z$ is small unless either $c_0=d_0=0$, in which case $Z$ contains the unprojection divisor $D_1=V(x,y,z,\nu)$, or $c_0=a_0 - d_0=b_0=0$, in which case $Z$ contains the unprojection divisor $D_2=V(x,y,z,\nu+a\eta)$. 

\paragraph{Isolated singularities.}
In case (i) $Q\in Y$ is a $\tfrac12(1,1,1)$ singularity and $Y$ must have isolated singularities. In case (ii) $Q\in Y$ is a $cA/2$ hyperquotient singularity. If both $a_0,c_0=0$ then $Z\supseteq\PP^1_{(\eta:\zeta)}$ and $Y$ becomes singular along this line.  

\paragraph{Description of $C\subset S_X$.} 
By explicitly resolving the $E_6$ Du Val singularity we see that, for case (i), the two numerical types of Table \ref{non_lci_table} are given by (i.a) $C_5+C_6$ if $d_0\neq0$ and (i.b) $C_1+C_4$ if $d_0=0$ and both $b_0,c_0\neq0$.

Similarly for case (ii) the three numerical types are given by (ii.a) $C_1+2C_5$ if $c_0\neq0$, (ii.b) $C_3+C_5$ if $c_0=0$ and all of $a_0,a_0-d_0,d_0\neq 0$ and (ii.c) $2C_4$ if $c_0=a_0-d_0=0$ and $b_0,d_0\neq 0$. 
%
%
%

In case (ii) a fourth numerical type is possible, given by $C_5 + 2C_6$ if $a_0=c_0=0$ and $d_0\neq 0$. The central fibre of $Y$ is small for this choice of $C$, so that $Y$ is a divisorial extraction, however $Y$ has non-isolated singularities.

\subsubsection{Type $\bm{E}_7$}

We assume that $C\subset S_X\subset X$ is of type $\bm{E}_7$. In this case $w\leq3$ and we split into the subcases: (i) $w=1$, (ii) $w=2$, (iii) $w=3$. 

\paragraph{The equations of $C$ and $Y$.} 
Our curve $C\subset X$ and $Y\subset X\times \PP(1,1,1,2)$ are given by the following equations:
\[ \bigwedge^2 \begin{pmatrix}
x & y^2 + z^3 & -(cy + dz^w) \\
y & x & ay + bz^w
\end{pmatrix} \quad\quad \Pf \begin{pmatrix}
\zeta & \nu & y\xi + c\eta & z^{3-w}\xi + d\eta \\
 & \xi & \nu + a\eta & b\eta \\
 & & -z^w & y \\
 & & & x
\end{pmatrix} \]
where we choose $w=w(D_0)$.

\paragraph{The general elephant $S_Y$.} 
In cases (i) and (ii) $\exc(\sigma |_{S_Y})=\PP^1_{(\xi:\zeta)}$ corresponding to $E_5$ in case (i) and to $E_1$ in case (ii). In case (iii) $\sigma|_{S_Y}$ is an isomorphism.

\paragraph{Small central fibre.}
In cases (i) and (ii) $Z$ is small unless either $b_0=d_0=0$ (contradicting $w=w(D_0)$) or $c_0=d_0=0$, in which case $Z$ contains an unprojection divisor $D_1=V(x,y,z,\nu)$. In case (iii) $Z$ is always small.

\paragraph{Isolated singularities.}
In case (i) $Q\in Y$ is a $\tfrac12(1,1,1)$ singularity and $Y$ must have isolated singularities. In case (ii) $Q\in Y$ is a hyperquotient singularity. If both $a_0,c_0=0$ then $Z\supseteq\PP^1_{(\eta:\zeta)}$ and $Y$ becomes singular along this line.  

\paragraph{Description of $C\subset S_X$.}
By explicitly resolving the $E_7$ Du Val singularity we see that the two possible numerical types for case (i) listed in Table \ref{non_lci_table} are given by (i.a) $C_1+C_6$ if $d_0\neq 0$ and (i.b) $C_4$ if $d_0=0$ and $b_0,c_0\neq 0$. The two types in case (ii) are given by (ii.a) $C_5+C_6$ if $c_0\neq 0$ and (ii.b) $C_2+C_6$ if $c_0=0$ and $a_0,d_0\neq 0$. Lastly, the two types in case (iii) are given by (iii.a) $3C_6$ if $c_0\neq 0$ and (iii.b) $2C_6+C_7$ if $c_0=0$ and $a_0\neq 0$. 
%
%
%
%
%
%
%
%
%
%

A third numerical type is possible in case (ii), given by $2C_1 + C_6$ if $a_0=c_0=0$ and $d_0\neq 0$. The central fibre of $Y$ is small in this case, so $Y$ is a divisorial extraction, however $Y$ has non-isolated singularities.

\subsection{Concluding remarks.}

In the $\bm{E}_6$ and $\bm{E}_7$ cases it is possible to take these calculations further and exhaust all of the possible cases by repeated serial unprojection. In both cases all Mori extractions exist in relative codimension $\leq5$ and the number of additional numerical types to consider is at most 12 for type $\bm{E}_6$ and at most 5 for type $\bm{E}_7$.

\end{document}